\documentclass[10pt]{article}
\title{\huge A Note on the Gauss Map of Complete Nonorientable Minimal 
Surfaces.}
\author{\Large Francisco J. L\'{o}pez  \thanks{Research partially
supported by
DGICYT grant number PB94-0796.} \\ Departamento de Geometr\'{\i}a y
Topolog\'{\i}a \\
Universidad de Granada \\18071 Granada, Spain \\ e-mail:{\tt 
fjlopez@goliat.ugr.es}
 \and \Large  Francisco Mart\'{\i}n $ ^{\ast}$ \\ Departamento de
 Geometr\'{\i}a y Topolog\'{\i}a \\
Universidad de Granada \\18071 Granada, Spain \\ e-mail:{\tt 
fmartin@goliat.ugr.es}}
\date{March 11, 1998}
\usepackage{latexsym}
\usepackage[latin1]{inputenc}

\newfont{\bb}{msbm10 at 10pt}
\newfont{\bt}{msbm10 at 6pt}
\newfont{\btt}{msbm10 at 5pt}
\newfont{\kt}{cmti10 at 6pt}
\newfont{\bg}{msbm10 at 15pt}

\def\r{\hbox{\bb R}}

\def\n{\hbox{\bb N}}
\def\c{\hbox{\bb C}}

\def\s{\hbox{\bb S}}
\def\d{\hbox{\bb D}}

\def\z{\hbox{\bb Z}}
\def\p{\hbox{\bb P}}

\newenvironment{proof}{\trivlist
\item[\hskip\labelsep{\em Proof}\,:]}{\hfill{$\Box$}\endtrivlist}

\def\dt{\hbox{\bt D}}

\newtheorem{lemma}{Lemma}

\newtheorem{theorem}{Theorem}

\newcommand{\df}{ \stackrel{\rm def}{=}}

\begin{document}
\maketitle
\begin{abstract}
We construct complete nonorientable minimal surfaces whose Gauss map omits two points
of $\r \p^2$. This result proves that Fujimoto's theorem is sharp in nonorientable
case.
\end{abstract}
\section{Introduction and Preliminaries} \label{sec:intro}
The study of the Gauss map of complete orientable minimal surfaces in $\r^3$ has achieved many important
advances and also has given rise to many problems in recent decades. The most interesting question is to determine
the size of the spherical image of such a surface under its Gauss map.

R. Osserman was the person who started the systematic development of this theory, and so,
in 1961 he proved that the set omitted by the image of a complete non flat orientable minimal surface by the
Gauss map has logarithmic capacity zero. In 1981 F. Xavier \cite{xavier}
proved that this set covers the sphere except six values at the most, and finally in 1988
H. Fujimoto \cite{fujimoto1,fujimoto2} obtained the best possible theorem,
and proved that the number of exceptional values of the Gauss map is four at the most.
An interesting extension of Fujimoto's theorem was proved in 1990 by  X. Mo and R. Osserman \cite{m-o}. They
showed that if the Gauss map of a complete orientable minimal surface takes on five
distinct values only a finite number of times, then the surface has finite total curvature.

There are many kinds of complete orientable minimal surfaces whose Gauss map omits four points of the
sphere. Among these examples we emphasize the classical
Sherk's doubly periodic surface and those described by K. Voss in \cite{voss}
(see also \cite{osserman}). The first author of this paper in \cite{lopez} constructs orientable examples
with non trivial topology.
 
Under the additional hypothesis of finite total curvature, R. Osserman \cite{osserman1} proved that the number of exceptional values
is three at the most. 

In the nonorientable case, the Gauss map of the two sheeted orientable covering surface  induces,
in a natural way, a {\em generalized Gauss map} from the nonorientable surface on the projective plane. From Fujimoto's theorem applied to the two sheeted orientable covering, this
generalized Gauss map omits two points of $\r \p^2$ at the most. 

It left open the following questions:
\begin{enumerate} \sl
\item Are there complete nonorientable minimal surfaces in $\r^3$ whose generalized Gauss
map omits two points of $\r \p^2$?
\item Are there complete non flat orientable minimal surfaces in $\r^3$ with finite total curvature whose Gauss map
omits three points of $\s^2$?
\item Are there complete nonorientable minimal surfaces in $\r^3$ with finite total curvature whose generalized Gauss
map omits one point of $\r \p^2$?
\end{enumerate}
Concerning the second problem, A. Weitsman and F. Xavier in \cite{w-x}  and Y.  Fang in \cite{fang}
obtained nonexistence results, provided that the absolute value of the total  curvature is less than
or equal to $16 \pi$ and $20 \pi$, respectively.

In this paper we give an affirmative answer to the first question, and prove:
\begin{quote}
{\bf Theorem} {\em There are complete nonorientable minimal surfaces in $\r^3$ whose generalized
Gauss map omits two points of the projective plane.}
\end{quote}

Our method of construction is somewhat explicit and very simple, and  it is based on a more  elaborate use of the Voss technique. 

Finally, we  briefly  summarize some of  the basic facts we will need in this paper.

Let $X:M \longrightarrow \r^3$ be a minimal immersion of a surface $M$ in three 
dimensional Euclidean space. Using 
isothermal parameters, $M$ has in a natural way  a conformal structure. When $M$ 
is orientable, we label $(g,\eta)$ as the Weierstrass data of $X$. Remember that 
the stereographic projection $g$ of the Gauss map of $X$ is a meromorphic function on $M$, and $\eta$ is a 
holomorphic 1-form on $M$.

Moreover,
$$ X= \mbox{Real}  \int ( \Phi_1,\Phi_2,\Phi_3),$$
where $\Phi_1=\frac{1}{2} \eta(1-g^2), \Phi_2=\frac{i}{2} \eta (1+g^2),  \Phi_3= \eta g$ are holomorphic 
1-forms on $M$ satisfying: $$\left( \sum^{3}_{j=1} |\Phi_{j}|^2 \right) (P) \neq 0, \quad \forall P \in M.$$
In particular, $\Phi_{j}$, $j=1,2,3,$ have no real periods on $M$. Furthermore, the 
Riemannian metric $ds^2$ induced by $X$ on $M$ is given by:
$$ds^2= \sum^{3}_{j=1} |\Phi_{j}|^2.$$
For more details see \cite{osserman}.

Consider now $X^{ \prime}:M^{ \prime } \longrightarrow \r^3$ a conformal minimal immersion of a 
nonorientable surface $M^{ \prime}$ in $\r^3$. Let $\pi_0:M \rightarrow M^{\prime}$, 
$I:M \rightarrow M$ denote the conformal oriented two sheeted covering of $M^{ \prime}$ 
and the antiholomorphic order two deck transformation for this covering, respectively.

If $(g,\eta)$ represents the Weierstrass data of $X=X^{ \prime} \circ \pi_0$, then it is not hard 
to deduce that :
\begin{equation} \label{eq:forms}
I^{\ast} (\Phi_{j})= \overline{\Phi_j}, \quad j=1,2,3.
\end{equation}
In particular, $g \circ I=I_0 \circ g$, where $I_0(z)=-1/\overline{z}$,
and so there is a unique  map $$G:M' \longrightarrow \r \p^2 \equiv 
\overline{\c}/\langle I_0 \rangle$$ satisfying 
$$G \circ \pi_0=g \circ p_0,$$
where $p_0:\overline{\c} \rightarrow \overline{\c}/\langle I_0 \rangle$ is the natural projection.
We call $G$ the {\em generalized Gauss map} of $X'$.

Conversely, given $(M,g,\eta)$ the Weierstrass representation of a minimal 
immersion $X$ of an orientable surface $M$ in $\r^3$, and given $I:M \rightarrow M$ an 
antiholomorphic involution without fixed points on $M$ satisfying (\ref{eq:forms}), then 
$X$ induces a minimal immersion $X^{ \prime}$ of $M^{ \prime}=M/ \langle I \rangle$ in $\r^3$ 
such that $X=X^{ \prime} \circ \pi_0$. For more details see \cite{meeks}.

Finally,  denote:
\begin{itemize}
\item $\d=\{ z\in \c \; :\; |z|<1 \}$,
\item $\d^\ast=\d-\{0\}$,
\item for each $R>1$, $A(R)=\{ z \in \c \; : \; 1/R<|z|<R \}$.
\end{itemize}

Throughout the proof of Theorem \ref{th:main}, we will use the following result:
\begin{theorem} \label{th:conmutativo}
Let $M$ be a Riemann surface with holomorphic universal covering space $\d$.
Then $M \cong \d $, $\d^\ast$, or $A(R)$, provided $\Pi_1(M)$ is commutative.
\end{theorem}
The proof of this theorem can be found in \cite[Chapter IV]{farkas-kra}.
\section{Main Theorem}
To obtain the result we have stated in the introduction, we need the following 
two Lemmas.

\begin{lemma} \label{lem:primero}
There exist $R>1$ and  holomorphic $1$-forms $\Phi_j$, $j=1,2,3$, on $A(R)$ such 
that:
\begin{enumerate}
\item $\Phi_1^2+\Phi_2^2+\Phi_3^2 \equiv 0.$
\item $|\Phi_1|^2+|\Phi_2|^2+|\Phi_3|^2 \neq 0.$
\item The metric $ds^2 \df  |\Phi_1|^2+|\Phi_2|^2+|\Phi_3|^2$ is complete.
\item The Gauss map $$g=-\frac{\Phi_1+i \Phi_2}{\Phi_3}$$ omits four points of 
the Riemann sphere
$\overline{\c}$.
\item $I^\ast (\Phi_j)=\overline{\Phi_j},$ $j=1,2,3,$ where $I:A(R) \rightarrow A(R)$ is given by
$I(z)=-1/\overline{z}$.
\end{enumerate}
\end{lemma}
\begin{proof}
Let $\alpha, \beta \in \c^\ast$, $\alpha \notin \{ \beta,-1/\overline{\beta} \}$, label 
$$M=\overline{\c} - \left\{\alpha,\beta,-\frac{1}{\, \overline{\alpha}
\,},-\frac{1}{\, \overline{\beta} \,} \right\}.$$
and consider the following Weierstrass representation on $M$:
\begin{equation} \label{eq:w1}
\widehat{g}=z, \; \; \widehat{\eta}=\frac{i dz}{(z-\alpha)(z-\beta) 
(\overline{\alpha} z+1)(\overline{\beta} z+1)}.
\end{equation}
If we define $\widehat{I}:M \rightarrow M$, $\widehat{I}(z)=-1/
\overline{z}$, then $\widehat{I}$ is an antiholomorphic involution
without fixed points, verifying:
\begin{equation} \label{eq:I1} \widehat{g} \circ \widehat{I}=-\frac{1}{\, \overline{\, \widehat{g} \, } \,}, \quad
\widehat{I}^\ast(\widehat{\eta})=-\overline{\widehat{\eta} \, \widehat{g}^2}. \end{equation}
Thus, if we define:
\begin{eqnarray*}
\widehat{\Phi}_1 & = & \frac{1}{2} (1-\widehat{g}^2) \widehat{\eta}, \\
\widehat{\Phi}_2 & = & \frac{i}{2} (1+\widehat{g}^2) \widehat{\eta}, \\
\widehat{\Phi}_3 & = & \widehat{g} \, \widehat{\eta}.
\end{eqnarray*}
then it is obvious, from (\ref{eq:I1}), that $ \widehat{I}^\ast( \widehat{\Phi}_j)=\overline{ \widehat{\Phi}_j}$.
Furthermore, these holomorphic 1-forms satisfy:
\begin{itemize}
\item $ \widehat{\Phi}_1^2+\widehat{\Phi}_2^2+\widehat{\Phi}_3^2 \equiv 0$,
\item $|\widehat{\Phi}_1|^2+|\widehat{\Phi}_2|^2+|\widehat{\Phi}_3|^2 \neq 0$,
\item The Riemannian metric $d \widehat{s}^2=|\widehat{\Phi}_1|^2+|\widehat{\Phi}_2|^2+|\widehat{\Phi}_3|^2$ is
complete in $M$.
\end{itemize}
On the other hand, the Uniformization Theorem  says us that the holomorphic
universal covering of $M$ is either $\c$ or the unit disc, $\d$ (see \cite[\S IV.4]{farkas-kra}). 
However, $\c$ is the conformal covering of only two non compact Riemann surfaces: $\c$ and $\c^\ast$
(see \cite[\S IV.6]{farkas-kra}). Thus,
the holomorphic universal covering of $M$ is $\d$.
 We label $\pi: \d \rightarrow M$ as the conformal covering map.

Let $\widetilde{I}$ be a lift of $\widehat{I}$ to $\d$, and denote $\widetilde{\Phi}_j = \pi^\ast(\widehat{\Phi}_j)$,
$j=1,2,3$.
It is clear that $\widetilde{I}^\ast(\widetilde{\Phi}_j)=\overline{\widetilde{\Phi}_j}$, $j=1,2,3$.

 Since $\widehat{I}$ is an antiholomorphic involution in $M$ without fixed points, then $\widetilde{I}^{2 k+1}$, $k \in \z$,
is an antiholomorphic transformation in $\d$ without fixed points too.

Let us see that  $\widetilde{I}^{2k}$, $k \in \z^\ast,$ has no fixed points in $\d$.
Indeed, note that  $\widetilde{I}^{2k}$, $k \in \z^\ast$, is a lift of the identity mapping in $M$.
Thus, if $\widetilde{I}^{2k}$ fixes
a point of $\d$, we infer that $\widetilde{I}^{2k}$ is the identity mapping $\mbox{\bf 1}_{\dt}$ in $\d$. 

Assume that there is $k>0$ such that $\widetilde{I}^{2 k}= \mbox{\bf 1}_{\dt}$.  Let 
$$k_0= \mbox{Minimum} \{k \in \n^\ast \;:\; \widetilde{I}^{2 k}= \mbox{\bf 1}_{\dt} \},$$
and observe that $k_0$ is the finite order of $\widetilde{I}^2$.
It is clear that $k_0>1$. Otherwise, $k_0=1$ and so  there would be 
antiholomorphic involutions without fixed points in $\d$, which is absurd.
Furthermore, from the definition of $k_0$, it is obvious that $\widetilde{I}^{2 k}$ has no fixed points, $0<k <k_0.$

Therefore, the quotient $\d / \langle \widetilde{I}^2 \rangle$ is a Riemann surface with fundamental group
isomorphic to $\z_{k_0}$. No such surface exists  (see for instance Theorem \ref{th:conmutativo}). 

This contradiction implies that $\widetilde{I}^{2 k}$, $k \in \z^\ast,$  has no fixed points and 
$ \langle \widetilde{I}^2 \rangle \cong \z$. In other words,  the map
$$ \zeta: \d \longrightarrow \d  / \langle \widetilde{I}^2 \rangle$$  is a
cyclic conformal covering and the  fundamental group of  $\d  / \langle \widetilde{I}^2 \rangle$ is
isomorphic to $\z$.

Using Theorem \ref{th:conmutativo} we deduce that $\d /\langle \widetilde{I}^2 \rangle$ is conformally equivalent to 
either $\d^\ast$ or $A(R)$, for a suitable $R>1$. 

The map $\widetilde{I}$ induces on
$\d / \langle \widetilde{I}^2 \rangle$ an antiholomorphic involution, $I$. Moreover, $\d / \langle \widetilde{I}^2 \rangle$
is in a natural way a covering of $M$, and $I$ is projected under this covering map on the original involution
$\widehat{I}$ on $M$. Since
$\widehat{I}$ has no fixed points in $M$, the same occurs for $I$ in $\d / \langle \widetilde{I}^2 \rangle$. 

However, any antiholomorphic involution in $\d^\ast$ extends to $\d$, and is the conjugate of a M\"{o}bius transformation
leaving $\d$ invariant and fixing $0$. In particular, any such map has infinitely many fixed points in $\d$. Hence,
we conclude that $\d / \langle \widetilde{I}^2 \rangle$ can not be conformally equivalent to $\d^\ast$, i.e., 
$\d / \langle \widetilde{I}^2 \rangle$  is conformally diffeomorphic to $A(R)$, for a suitable  $R>1$.

If we look at  $I$ as an antiholomorphic involution in  $A(R)$, then elementary arguments of complex analysis give that
$ I(z)=-1/\overline{z}$, $\forall z \in A(R)$.

On the other hand, as $(\widetilde{I}^2)^\ast (\widetilde{\Phi}_j)=\widetilde{\Phi}_j$, then $\widetilde{\Phi}_j$
can be induced in the quotient $\d / \langle \widetilde{I}^2 \rangle$,  $j=1,2,3$. The corresponding holomorphic 1-forms
on $\d / \langle \widetilde{I}^2 \rangle$ are denoted as $\Phi_1$, $\Phi_2$, and $\Phi_3$, and they obviously satisfy
{\em 1, 2, 3} and {\em 5} in the lemma statement.

Finally, the meromorphic function $$g=-\frac{\Phi_1+i \Phi_2}{\Phi_3},$$ clearly omits the points
$\alpha$, $\beta$, $-1/\overline{\alpha}$, and $-1/\overline{\beta}$, and {\em 4} holds. This concludes
the proof.

\end{proof}

\begin{lemma} \label{lem:segundo}
There exists a rational function $f:\overline{\c} \rightarrow \overline{\c}$ satisfying:
\begin{enumerate}
\item The only poles of $f$ are $0$ and $\infty.$
\item $f \circ I=\overline{f}.$
\item $f(z) \neq 0$, provided that $|z|=1.$
\item $\displaystyle \mbox{Residue} \left( \frac{f(z)}{z} dz, 0 \right) =0.$
\end{enumerate}
\end{lemma}
\begin{proof}
Define $$f:\c \longrightarrow \c,$$ $$f(z)=\frac{(z-m_1) (z-m_2) (m_1 z+1) (m_2 z+1)}{z^2},$$
where $m_1, m_2 \in \r.$

We have $$  \mbox{Residue} \left( \frac{f(z)}{z} dz, 0 \right) =   (1-m_1^2)(1-m_2^2) -2 m_1 m_2.$$

The choice $m_1=2$ and $m_2=\frac{2 +\sqrt{13}}{3} $ completes the proof.
\end{proof}
Now we are able to prove the main result of this paper.

\begin{theorem} \label{th:main}
There exist complete nonorientable minimal surfaces in $\r^3$ whose
generalized Gauss map omits  two points of $\r \p^2$.
\end{theorem}
\begin{proof}
Take $A(R)$, $\Phi_1$, $\Phi_2$, and $\Phi_3$ as in Lemma \ref{lem:primero}, and 
$f$ as in Lemma \ref{lem:segundo}. Put
$$ \Phi_j=  \varphi_j(z) \frac{dz}{z}, $$
and write
$$\varphi_j(z)= a_{j \, 0}+ \sum_{n > 0} \left(a_{j \, n} \, z^n+(-1)^{n+1} \; \overline{a_{j \, n}} \, z^{-n} \right) , \quad
a_{j \, 0} \in i \r,$$
the Laurent  series expansion of $\varphi_j$, $j=1,2,3$.

Observe that 
$$f(z)=\sum_{n=1}^m \left( b_n \, z^n+ (-1)^n \, \overline{b_n} \, z^{-n} \right),$$
where $m \in \n^\ast.$
Let $k \in \n$, $k$ odd, $k>m$, and notice that: 

\begin{equation} \label{eq:uno}
\mbox{Residue} \left( \left[ \sum_{n > 0} \left(a_{j \, n} \, z^{k n}+ (-1)^{n+1} \,  \overline{a_{j \, n}} \, z^{-k n} \right) \right] \,
f(z) \frac{dz}{z},0  \right)=0, \;\; j=1,2,3.
\end{equation}

Furthermore, it is obvious from Lemma \ref{lem:segundo} 
\begin{equation} \label{eq:dos}
\mbox{Residue} \left( a_{j \, 0} \, f(z) \, \frac{dz}{z},0 \right) =0, \;\; j=1,2,3. 
\end{equation}

Consider the covering  $T_k:A(\sqrt[k]{R}) \rightarrow A(R)$, $T_k(z)=z^k,$ and define the 
holomorphic $1$-forms on $A(\sqrt[k]{R})$:
$$\Psi_j= f(z)  \, T_k^\ast (\Phi_j)=k \,  f(z) \,  \varphi_j(z^k) \,  \frac{dz}{z}, \; \; j=1,2,3.$$
Taking into account (\ref{eq:uno}) and (\ref{eq:dos}), we deduce that $\Psi_j$ is exact, 
$j=1,2,3.$ 

Moreover, it is clear that:
$$\sum_{j=1}^3 \Psi_j^2 \equiv 0,$$ and since $k$ is odd,
\begin{equation}
I^\ast(\Psi_1,\Psi_2,\Psi_3)  =  \left( \overline{\Psi_1}, \overline{\Psi_2}, \overline{\Psi_3} \right), \label{eq:psi2}
\end{equation}
where $I:A(\sqrt[k]{R}) \rightarrow A(\sqrt[k]{R})$  is the lift of the former involution in $A(R)$, that
keeps being the map $I(z)=-1/\overline{z}.$

Note that $\lim_{k \to \infty} \sqrt[k]{R}=1$, and remember that the zeroes of $f$ are not in $\s^1$.
Then, taking $k$ large enough, we can guarantee that $f$ never vanishes in the closure of $A(\sqrt[k]{R}).$
So, as the only poles of $f$ are $0$ and $\infty$, there exist $c>1$ such that
$$ \frac{1}{c} < |f(z)|<c, \quad \forall z \in A(\sqrt[k]{R}). $$
Therefore, $\sum_{j=1}^3 | \Psi_j|^2  \neq  0,$
and if we define $ds_0^2=|\Psi_1|^2+|\Psi_2|^2+|\Psi_3|^2$, one has:
$$\frac{1}{c^2} \, T_k^\ast( ds^2) \leq ds_0^2 \leq c^2 \, T_k^\ast (ds^2).$$
Since $ds^2$ is complete,  the same occurs for the metrics $T_k^\ast (ds^2)$ and  $ds_0^2$. 

Summarizing, the minimal immersion  $$X:A(\sqrt[k]{R}) \longrightarrow \r^3,$$
$$X(z)= \mbox{Real} \left( \int_1^z (\Psi_1,\Psi_2,\Psi_3) \right), $$
is well defined, complete, and its Gauss map $g \circ T_k$ omits four
points of $\overline{\c}$.

From (\ref{eq:psi2}), $X$ induces a minimal immersion of the M\"{o}bius strip
$A(\sqrt[k]{R})/\langle I \rangle $ in $\r^3$, and so the Theorem is proved.

\end{proof}

\end{document}